\documentclass[11pt]{article} 
\usepackage{amsmath,amsthm,amscd,amssymb}
\usepackage{latexsym}
\usepackage{epsf}
\usepackage{epsfig}
\newcommand{\diff}{\operatorname{Diff}}

\newcommand{\vol}{\operatorname{Vol}}
\renewcommand{\div}{\operatorname{div}}
\newcommand{\R}{\mathbb{R}}

\newcommand{\N}{\mathbb{N}}

\newcommand{\w}{\omega}
\newcommand{\om}{\omega}

\newcommand{\ep}{\varepsilon}

\newcommand{\LM}{L(M)}
\newcommand{\LN}{L(N)}
\newcommand{\SA}{\mathcal{A}}
\newcommand{\SQ}{\mathcal{Q}}

\newcommand{\SM}{\mathcal{M}}
\newcommand{\SP}{\mathcal{P}}

\newcommand{\II}{\operatorname{II}}

\theoremstyle{plain}
\newtheorem{theorem}{Theorem}[section]
\newtheorem{lemma}[theorem]{Lemma}

\newtheorem{proposition}[theorem]{Proposition}

\theoremstyle{definition}
\newtheorem{definition}[theorem]{Definition}

\theoremstyle{remark}

\newtheorem{case[theorem]}{Case}

\title{Minimal Distortion Bending and Morphing of Compact Manifolds 
\footnote{This research is supported by the NSF 
Grant NSF/DMS-0604331}}
\author{Oksana Bihun\footnote{ Corresponding Author E-mail: oksana@math.missouri.edu}\hspace*{.05in} and Carmen Chicone\footnote{E-mail: carmen@math.missouri.edu}\\Department of Mathematics\\University of
Missouri-Columbia\\Columbia, Missouri 65211, USA}

\begin{document}
\maketitle
%\date{\today}
geometric optimization, minimal morphing

\begin{abstract}
 Let $M$ and $N$ be
compact smooth oriented Riemannian $n$-manifolds  without boundary 
embedded in $\mathbb{R}^{n+1}$. Several problems about minimal distortion 
bending and morphing of $M$ to $N$ are posed.
Cost functionals that measure 
distortion due to stretching or bending produced by a diffeomorphism 
$h:M \to N$  are defined, and new results on the existence of minima of these cost functionals 
are presented. In addition, the definition of a morph between two manifolds $M$ and $N$ is
 given, and the theory
of minimal distortion morphing of compact manifolds is reviewed.
\end{abstract}

\section{Introduction}
Two diffeomorphic compact embedded hypersurfaces admit infinitely many 
diffeomorphisms between them, which we view as prescriptions for bending one hypersurface 
into the other. We ask which diffeomorphic bendings have minimal distortion 
with respect to some natural bending energy functionals. More precisely,
let $M$ and  $N$ be diffeomorphic compact and connected smooth oriented 
$n$-manifolds without boundary
embedded in $\R^{n+1}$. The manifolds $M$ and $N$ inherit Riemannian metrics $g_M$ and $g_N$ and corresponding volume forms  from the usual metric and orientation of $\R^{n+1}$. Although we  will use this structure here, only the
existence of the metrics on $M$ and $N$ is essential.  We pose the problem of
bending $M$ into  $N$ via a diffeomorphism 
$h:M \to N$ so that the distortion produced by $h$ is minimal with respect to some functional that measures bending or stretching (cf. the problem of optimal development of surfaces~\cite{YPM, LYS}).

The problem of minimal distortion bending of manifolds may be considered a special case of the problem of minimal morphing. A morph is defined to be a transformation between 
two shapes through a set of intermediate shapes.  A minimal morph is such a 
transformation that minimizes distortion.  
There are important applications of minimal morphing in  
manufacturing~\cite{LYS, YPM},
computer graphics~\cite{W,W1}, movie making~\cite{HLZ}, 
and mesh construction~\cite{HH,LDSS}. 
We will formulate and solve a problem about 
the existence of minimal morphs with respect to stretching for $n$-dimensional manifolds.

\section{Distortion Minimal Bending}
Let $(M, g_M)$ and $(N, g_N)$ be smooth compact Riemannian $n$-manifolds without
boundary with all the additional properties stated in the introduction. 
By $\vol(M)$ we denote the volume of $M$. Also, let $\diff(M,N)$ denote the set 
of all  diffeomorphisms between $M$ and $N$. 

Our first functional  measures distortion due to stretching. For a point $p \in M$, we define the distortion due to stretching at
$p$ as the infinitesimal relative change of volume. More precisely, let
$\{A_k\}_{k=1}^\infty
\subset M$ be a sequence of open neighborhoods of the 
point $p$ that shrink to the point as $k\to\infty$. For example, one can choose $A_k=B(p,\frac{1}{k})\cap M$, where 
$B(p,R)$ is the open ball of radius $R$ in $\R^{n+1}$ centered at $p \in M
\subset \R^{n+1}$.
\begin{definition}
The \emph{distortion due to stretching} produced by a diffeomorphism $h \in
\diff(M,N)$ at a point $p\in M$ is defined to be
$$\xi(p)=\lim \limits_{k \to \infty} 
\frac{\big|\int_{h(A_k)} \w_N\big|-\big|\int_{A_k} \w_M \big| }
{\big| \int_{A_k} \w_M \big|}=\big|J(h)(p)\big|-1,$$
where $J(h)$ is the Jacobian of $h$ with respect to the (Riemannian) volume
forms $\w_M$ and $\w_N$.
The functional $\Phi_1:\diff(M,N)\to \R_+$
is defined by 
\begin{equation}
\label{defPhi1}
\Phi_1(h)=\int_M\Big(\big|J(h)\big|-1\Big)^2\,  \w_M.
\end{equation} 
\end{definition}

The following results are proved in~\cite{RMJpaper}.
\begin{lemma}
\label{le:Jconst}
A diffeomorphism $h \in \diff(M,N)$ is a critical point of $\Phi_1$ if and 
only if $J(h)(m)=\frac{\vol(N)}{\vol(M)}$ for all $m \in M$.
\end{lemma}
\begin{theorem}[Existence of minimizers for $\Phi_1$]
\label{th:minphi}
If $(M, g_M)$ and $(N, g_N)$  are diffeomorphic 
compact 
connected oriented Riemannian $n$-manifolds without boundary,  then
there exists a minimizer of the functional $\Phi_1$ over the class 
$\diff (M,N)$ and
 the minimum value of $\Phi_1$ is
\begin{equation*}
\Phi_1^{\min}=\frac{\big( \vol(M)-\vol(N)\big)^2}{\vol(M)}.
\end{equation*}
\end{theorem}

The functional $\Phi_1$ is invariant with respect to compositions with volume
preserving maps: $\Phi_1(h \circ k)=\Phi_1(h)$ provided that $k \in \diff(M)$ is
volume preserving (it has the Jacobian $J(k)=1$). Therefore, the minimizer of
$\Phi_1$ is not unique. 

For a vector bundle $V$ over $M$, we denote by $\Gamma(V)$ the space of all 
sections of $V$;
$\Gamma^r(V)$ denotes the space of all $C^r$ sections of $V$.
Let $T^{(0,2)}(M)$ denote the vector bundle of covariant order-two tensors over 
$M$  (see \cite{AMR}).
In order to measure distortion with respect to bending, we introduce the 
strain tensor field $S:\diff(M,N) \to \Gamma(T^{(0,2)}(M))$ by
$S(h)=h^\ast g_N-g_M$,  where $h^\ast g_N$ is the pull-back of the metric $g_N$ 
by $h$. 

\begin{definition}
The \emph{deformation energy functional} $\Phi_2:\diff(M,N)\to \R_+$ is  given by 
$$
\Phi_2(h)=\int_M \|h^\ast g_N-g_M\|^2\w_M,
$$
where the fiber norm $\|\cdot\|$ on the bundle $T^{(0,2)}(M)$ 
is induced by the fiber metric $G:=g_M^\ast \otimes g_M^\ast$ 
(see \cite{KN}).
\end{definition}

Let $\Gamma^\infty(TM)$ denote the space of
$C^\infty$ sections of the tangent bundle of $M$.
Fix a diffeomorphism $h \in \diff(M,N)$. In order to derive the Euler-Lagrange 
equation, 
we consider all variations 
of the diffeomorphism $h$ of the form $h\circ\phi_t $, where $\phi_t$ is the
flow of a vector field $Y \in \Gamma^\infty(TM)$. Because the tangent space $T_h
\diff(M,N)$ can be identified with the space of all sections $\Gamma^\infty (h^{-1}TN)$
of the pull back bundle $h^{-1}TN$ over $M$, every smooth variation of the
diffeomorphism $h$ can be represented in the form $h \circ \phi_t$ (see~
\cite{OpusculaPaper} for a more detailed description).

The diffeomorphism $h \in \diff(M,N)$ is a critical point of $\Phi_2$ if
\begin{equation}
\label{EL3}
\frac{d}{dt} \Phi_2(h \circ \phi_t)|_{t=0 }=
D\Phi_2(h) h_\ast Y=2 \int_M G(h^\ast g_N-g_M, L_Y h^\ast g_N)=0
\end{equation}
for all $Y\in \Gamma^\infty(TM)$, where $h_\ast Y$ is the push forward of $Y$ by
$h$, and $L_Y$ is the Lie derivative in the direction of $Y$ (see
\cite{AMR}).

Let $\nabla$ be the Riemannian connection on $M$ generated by the Riemannian
metric $g_M$ with Christoffel symbols
$\Gamma_{ij}^k$ (see~\cite{Hicks}).  The connection $\bar{\nabla}$ with Christoffel symbols
$\bar{\Gamma}_{ij}^k$ is the Riemannian connection of the metric $h^\ast g_N$ on $M$.

\begin{definition}
\label{defB}
Define $B(h)=(h^\ast g_N- g_M)^\#$. That is, for each $p \in M$, the tensor
$B(h)(p)$  of type $(2,0)$ 
is defined as the tensor $(h^\ast g_N- g_M)(p)$ with its indices raised. 
\end{definition}

\begin{definition}
\label{defS}
Define the bilinear form $A(h)$ to be
\begin{equation}
\label{SS}
A(h)(X,Y)=\bar{\nabla}_X Y-\nabla_X Y
\end{equation}
for $X,Y \in \Gamma^\infty(TM)$ (see \cite{KN}).
\end{definition}
It is easy to prove that $A(h)$ is a tensor field of type $(1,2)$ on $M$ with components
$$A(h)^m{}_{kp}=\bar{\Gamma}_{kp}^m -\Gamma_{kp}^m.$$

\begin{lemma}
\label{fundlemma1}
The first variation of the functional $\Phi_2$ in the direction $Y \in \Gamma^{\infty}(TM)$ is
given by
\begin{equation}
\label{DPhiY1}
D\Phi_2(h)(h_\ast Y)=-4\int_M g_M(\div B(h)+ A(h):B(h),Y)\,\w_M,
\end{equation}
 where $B(h):A(h)$ is the contraction of the
tensor fields $B(h)$ and $A(h)$ (see \cite{MH}).
Moreover, 
$h$ is a critical point of the functional $\Phi_2$ if and only if 
\begin{equation}
\label{ELBump}
\div B(h)+A(h):B(h)=0.
\end{equation}
\end{lemma}
\noindent The Euler-Lagrange equation for the functional $\Phi_2$ is the system of nonlinear
partial differential equations \eqref{ELBump}. 

Let $h_R: \mathbb{R}^{n+1}
\to  \mathbb{R}^{n+1}$ be the radial map given by 
$h_R(x)=R x$ for some number $R>0$ and for all $x \in \R^{n+1}$. It is easy to check that if $N=RM$ is a
rescaling of the manifold $M$, then the map $h=h_R\circ f$ satisfies the Euler
Lagrange equation~\eqref{ELBump}, whenever $f \in \diff(M)$ is an isometry on $M$.

The following results on minimizing $\Phi_2$ in the one-dimensional case are proved in \cite{OpusculaPaper}.
\begin{proposition}
\label{exist1dim}
(i) Suppose that $M$ and $N$ are smooth simple closed curves in $\R^2$ with arc lengths $\LM$ and $\LN$ and base points $p\in M$ and $q \in N$; 
 $\gamma$ and $\xi$ are the  arc length parametrizations of $M$ and $N$ with $\gamma(0)=p$ and $\xi(0)=q$ that induce positive orientations; 
and, the functions $v$ and $w$ are defined by
$v(t)=\LN/\LM t$  and $w(t)=-\LN/\LM t+\LN$ for all $t \in [0,
\LM].$
If $\LN \geq\LM$, then
the functional $\Phi_2$ has exactly two
minimizers in the admissible set $$\SA=\{h \in \diff(M,N):h(p)=q\}:$$
the orientation preserving minimizer 
\[
h_1=\xi \circ v \circ \gamma^{-1}
\]
and the orientation reversing minimizer
\[
h_2=\xi \circ w \circ \gamma^{-1}
\]
(where we consider $\gamma$ as a function defined on  
$\big[0,\LM \big)$ 
so that $\gamma^{-1}(p)=0$). Moreover, 
the minimum value of the functional $\Phi_2$ is
\begin{equation}
\Phi_2^{\min}=\frac{(\LN^2-\LM^2)^2}{\LM^3}.
\end{equation}
\end{proposition}

\begin{proposition}
\label{cor1} Assume the notation of the  proposition~\ref{exist1dim}.\\
(i) If $\LN<\LM$, then the functional $\Phi_2$ has no minimum in the admissible set
\[
\SQ=\{h \in C^{\infty}(M,N): \mbox{$h$ is orientation preserving
and $h(p)=q$} \}.
\]
(ii) If 
$\frac{\LN}{\LM}<\frac{1}{\sqrt{3}}$,  then the functional $\Phi_2$ has 
no minimum in the admissible set $\SA=\{h \in \diff(M,N):h(p)=q\}$.
\end{proposition}

The main ingredients for a proof of (i) in proposition~\ref{cor1} are simply illustrated.  The curve $M$ is  wrapped around the curve $N$ without stretching  and  the excess is removed. This wrapping function can be expressed in the form 
$h=\xi\circ u\circ \gamma^{-1}$ (in the notation of
proposition~\ref{exist1dim}), where $u:[0,\LM] \to \R$ is a discontinuous
piecewise linear function. The
function $h$ is not smooth; but, it is possible to approximate it by a 
minimizing sequence $\{h^k\}_{k=1}^\infty \subset \SQ$ whose deformation 
energies $\Phi_2(h^k)$ converge to $\Phi_2(h)=0$. On the other hand,
$\Phi_2(f)>0$ for all $f \in Q$.
The proof of (ii) uses the second variation of $\Phi_2$.

By proposition~\ref{cor1}, we see that even in the one-dimensional
case  the functional $\Phi_2$ exhibits nontrivial behavior: the minimum does
not always exist, and the existence depends on properties of the curves $M$ and
$N$.

The general problem of the existence of minimizers for the functional $\Phi_2$ is open. On the other hand, we have solved the problem for the case where
$M$ and $N$ are Riemann spheres or compact Riemann surfaces of genus greater
than one. Let $\mathcal{H}(M,N)=\{h \in \diff(M,N): h \mbox{ is a holomorphic map}\}$.
\begin{theorem}
\label{RiemSph}
(i) Let $h_R: \mathbb{R}^3\to  \mathbb{R}^3$ be the radial map given by 
$h_R(p)=R p$ for some number $R>0$. If  $M=S^2 \subset \R^3$ and  $N=h_R(M)$,  
then  $h:=f \circ h_R$ is a global minimum 
of the functional $\Phi_2$, restricted to the admissible set 
$\mathcal{H}(M,N)$,  whenever $f$ is 
an isometry of $N$.

(ii) Let $M$  and $N$ be compact Riemann surfaces. If
 $\mathcal{H}(M,N)$ is not empty and
the genus of $M$ is at least two, then
there exists a minimizer of the functional $\Phi_2$ in $\mathcal{H}(M,N)$. 
\end{theorem}

The general problem of minimization of the functional $\Phi_2$ seems 
to be very difficult because the admissible set
is an infinite-dimensional manifold $\diff(M,N)$ whose structure is not completely understood. In 
theorem~\ref{RiemSph}, the admissible set is a finite-dimensional homogeneous
space in case $M$ and $N$ are two-spheres and a finite group in case $M$ is a compact Riemannian manifold of genus greater than one.  A natural idea is to  
reformulate the problem of minimal distortion bending in such a way that 
the admissible set is a linear space.

Fix a diffeomorphism $f \in \diff(M,N)$. Every diffeomorphism $h:M \to N$ can be represented in the form $h=f \circ \phi$, where $\phi \in \diff(M)$. For simplicity of notation, let 
$g_2:= f^\ast g_N$ and $g_1:=g_M$. To measure
the deformation produced by $h=f\circ \phi \in \diff(M,N)$, we use the 
strain tensor field
$S(\phi)=\phi^{\ast} g_2 - g_1$. In other words, the problem reduces 
to minimization of the deformation energy produced by some class of diffeomorphisms 
$\phi:(M,g_1) \to (M,g_2)$ in $\diff(M)$.

The tangent bundle $TM$ is equipped with the Riemannian metric $g_1$. 
 Let $W^{k,2}(TM)$ be the $(k,2)$-Sobolev space of sections of the
 tangent bundle
$TM$ (see \cite{Wehrheim}). We choose the number of generalized derivatives 
 $k \in \N$ large enough so that the Sobolev space 
$W^{k,2}(TM)$ is embedded into  the space $\Gamma^2(TM)$ of all $C^2$ sections of $TM$
and some additional estimates hold.
Consider the space $H=L^{2}\big([0,1];W^{k,2}(TM)\big)$ of time dependent
vector fields $v:M\times[0,1]\to \Gamma(TM)$. The space $H$ is a Hilbert 
space equipped with the norm
$$\langle v,w \rangle_H=\int_0^1\langle v(\cdot, t), w(\cdot,t)
\rangle_{W^{k,2}(TM)} \,dt.$$  

Every vector field $v \in H$ generates a diffeomorphism on $M$ in the following
sense. The nonautonomous ordinary differential equation
\begin{equation}
\label{IVP2}
\frac{d q}{dt}=v(q,t),
\end{equation}
is solved (on the compact manifold $M$) by an evolution operator 
$\eta^v(t; s ,p)$ that satisfies the Chapman-Kolmogorov conditions 
(see~\cite{Chicone})
and is such that $\eta^v(s; s ,p)=p$ for every $p \in M$. 
The function $\phi^v:M\to M$ given by
$\phi^v(p)=\eta^v(1;0,p)$  is called the time-one map of the evolution
operator $\eta^v$; it is a
diffeomorphism on the manifold $M$. 

We define the distortion
energy functional $E:H \to \R_+$ to be
\begin{equation}
\label{DefE}
E(v)=\|v\|_H^2+\int_M\|(\phi^v)^\ast g_2- g_1\|^2 \w_M+
\int_M\|(\phi^v)^\ast \II_2- \II_1\|^2 \w_M,
\end{equation}
where $\II_i$ is the second fundamental form on $M$ associated with $g_i$,
$i=1,2$ (see~\cite{Hicks}).
This functional incorporates strain (which is intrinsic to the manifold $M$) and
bending (which is extrinsic).

\begin{theorem}[Existence of minimizers for $E$]
There exists a minimum of the functional $E$ in the space $H$.
\end{theorem}
The proof uses the direct method of the calculus of variations as well as convergence
properties of the evolution operators $\eta^{v^{l}}(t;s,x)$ generated by weakly convergent
sequences $\{v^l\}_{l=1}^\infty$ of time dependent vector fields in $H$.

\section{Distortion Minimal Morphing}
\begin{definition}
\label{df:morph}
Let $M$ and $N$ be compact connected oriented $n$-dimensional smooth 
manifolds without boundary embedded in
$\R^{n+1}$. A $C^1$ function
$F:[0,1] \times M \to \R^{n+1}$ is a morph from $M$ to $N$ if the
following conditions hold:
\begin{itemize}
\item[(i)] $p \mapsto F(t,p)$ is a diffeomorphism onto its image for each $t
\in I=[0,1]$;
\item[(ii)]  the image $M^t=F(t,M)$ is an $n$-dimensional manifold possessing all the
properties of $M$ and $N$ mentioned above;
\item[(iii)]  $p \mapsto F(0,p)$ is a diffeomorphism of $M$;
\item[(iv)]  the image of the map $p \mapsto F(1,p)$ is $N$.
\end{itemize}
We denote the set of all $C^2$ morphs between the manifolds $M$ and $N$ by $\SM(M,N)$.
\end{definition}

For simplicity, we will only consider morphs $F$ such that $p \mapsto F(0,p)$ is the identity map. We assume that each manifold $M^t=F(t,M)$ (with  $M^0=M$ and $M^1=N$) is equipped with the volume form $\om_t=i_{\eta_t}\Omega$,
where $\Omega=dx_1 \wedge dx_2 \wedge \ldots \wedge dx_{n+1}$ is the standard
volume form on $\R^{n+1}$ and $\eta_t:M^t \to \R^{n+1}$ is the outer unit normal
vector field on $M^t$ with respect to the usual metric on $\R^{n+1}$. Also,  as
a convenient notation, we use $f^t=F(t,\cdot):M \to M^t$.
\begin{definition}
The
functional
$\Phi_1^{s,t}:\diff(M^s,M^t) \to \R_+$ is defined by formula \eqref{defPhi1},
were $M$ and $N$ are replaced with $M^s$ and $M^t$ respectively.
A morph $F$ is \emph{distortion pairwise minimal}
 (or, for brevity, \emph{pairwise minimal}) if $f^{s,t}=f^t \circ (f^s)^{-1}:M^s
\to M^t$ minimizes the functional $\Phi_1^{s,t}$ for every  $s, t \in[0,1]$. We denote the set of all 
$C^2$ distortion
pairwise minimal morphs between manifolds $M$ and $N$ by $\SP \SM(M,N)$. 
\end{definition}

Using lemma~\ref{le:Jconst}, it is easy to derive a necessary and sufficient condition
for pairwise minimality.
\begin{proposition}
\label{pr:pairwise}
Let $M=M^0$ and $N=M^1$ be $n$-dimensional manifolds as in definition~\ref{df:morph}
equipped with the (respective) volume forms $\om_0$ and $\om_1$. A morph $F$ between $M$ and $N$ is distortion pairwise minimal 
if and only if
\begin{equation}
\label{ConstPairwise}
\frac{J(f^t)(m)}{\vol(M^t)}=\frac{1}{\vol(M)}
\end{equation}
for all $t \in [0,1]$ and $m\in M$, where $J(f^t)$ is the Jacobian of $f^t$ with
respect to the volume forms $\w_0$ and $\w_t$.
\end{proposition}

The following proposition states the existence of pairwise minimal morphs. It can be proved by rescaling morphs between  $M$ and $N$, which are not
necessarily pairwise minimal, to conform to property \eqref{ConstPairwise}. Moser's theorem on volume forms (see \cite{M}) plays a
crucial role in the proof.
\begin{proposition}
\label{pr:pwminimal}
Let $M$ and $N$ be $n$-dimensional manifolds as in definition~\ref{df:morph}. 
If $M$ and $N$ are connected by a $C^2$ morph, then there is a distortion pairwise
minimal morph between them.
\end{proposition}

Having the preliminary study of pairwise minimal morphs at hand, we define
minimal morphs.
\begin{definition}
The \emph{infinitesimal distortion} of a $C^2$ morph $F$ from $M$ to $N$ at 
$t\in [0,1]$ is 
$$
\ep^F(t)= \lim \limits_{s \to t} \frac{E^{s,t}}{(s-t)^2}=
\int_M \frac{\big(\frac{d}{dt} J(f^t) \big)^2}{J(f^t)} \om_M,
$$
where $E^{s,t}=\Phi_1^{s,t}(f^{s,t})$ is the distortion energy of the transition 
map $f^{s,t}$.
The \emph{total distortion functional} $\Psi$ defined on such morphs is given by \begin{equation}
\label{eq4:1}
\Psi(F)=\int_0^1 \ep^F(t) \,dt= \int_0^1 
\Big(\int_M \frac{\big(\frac{d}{dt} J(f^t) \big)^2}{J(f^t)}
\om_M\Big)dt.
\end{equation}
\end{definition}

The following proposition implies that it suffices to minimize the functional
$\Psi$  over the class $\SP\SM(M,N)$ of all pairwise minimal morphs instead of
the class $\SM(M,N)$ of all $C^2$ morphs.
\begin{proposition}
\label{Cor:PairwInf}
\begin{itemize}
\item[(i)] The following inequality holds:
\begin{equation}
\inf_{G \in \SP\SM(M,N)} \Psi(G) \leq \inf_{P \in \SM(M,N)} \Psi(P).
\end{equation}
\item[(ii)]
If there exists a minimum $F$ of the total distortion functional $\Psi$ over 
the class $\SP\SM(M,N)$, then $F$ minimizes the functional $\Psi$ over the class $\SM(M,N)$ as well; in fact,
\begin{equation}
\Psi(F)=\min_{G \in \SP\SM(M,N)} \Psi(G) =\min_{P \in \SM(M,N)} \Psi(P).
\end{equation}
\end{itemize}
\end{proposition}

Using proposition~\ref{pr:pairwise}, it is easy to recast the functional $\Psi$
into a simpler form:
\begin{lemma} 
\label{Lemma:pairwEn}
The total distortion of a $C^2$ pairwise minimal morph
$F$ from $M$ to $N$ is
\begin{equation}
\label{eq4:2}
\Psi(F)=\int_0^1 \frac{\big(\frac{d}{dt}\vol(M^t)\big)^2}{\vol(M^t)}\, dt.
\end{equation}
\end{lemma}
The latter form of the functional $\Psi$ and proposition~\ref{Cor:PairwInf}
allow us to solve the problem of minimization of the total distortion functional
$\Psi$ over the class of all $C^2$ morphs $\SM(M,N)$. In order to solve the
problem, we
minimize the auxiliary functional
\begin{equation}
\label{auxfct}
\Xi(\phi)=\int_0^1 \frac{\dot{\phi}^2}{\phi}\, dt
\end{equation}
over the admissible set $$Q=\big\{\phi \in C^1\big([0,1];\R_+\big):
\phi(0)=\vol(M), \phi(1)=\vol(N)\big\}.$$

The following theorem---our main result on distortion minimal 
morphing---is proved using  
proposition~\ref{Cor:PairwInf} and lemma~\ref{Lemma:pairwEn}.

\begin{theorem}
Let $M$ and $N$ be two $n$-dimensional manifolds satisfying the assumptions 
of definition~\ref{df:morph}. If $M$ and $N$ are connected by a $C^2$ morph,  
 then they are connected by a minimal morph. The minimal value of
$\Psi$ is
\begin{equation}
\label{PhiMin11}
\min_{F \in \SM(M,N)} \Psi(F) = 4
\big(\sqrt{\vol(N)}-\sqrt{\vol(M)}\big)^2.
\end{equation}

\end{theorem}


\begin{thebibliography}{xxxxx}
\bibitem{AMR} R. Abraham, J. Marsden and T. Ratiu, \emph{Manifolds, Tensor
Analysis, and Applications}, Springer-Verlag New York Inc., 1988.
\bibitem{OpusculaPaper} O. Bihun and C. Chicone, {Deformation Minimal Bending of Compact 
Manifolds: Case of Simple Closed Curves}, \emph{Opuscula Mathematica}, to appear in January 2008,\\ 
http://arxiv.org/abs/math.OC/0701901.
\bibitem{RMJpaper} O. Bihun and C. Chicone, 
{Distortion Minimal Morphing I: The Theory For Stretching}, Preprint 2006, \\
http://arxiv.org/abs/math.DG/0605668.
\bibitem{Chicone} C. Chicone, \emph{Ordinary Differential 
Equations with Applications}, Texts in Applied Mathematics, New York: 
Springer-Verlag, 2006.
\bibitem{HH}H. Hoppe,  Progressive meshes.
\emph{ACM SIGGRAPH} (1996) 99--108.
\bibitem{LDSS} A. Lee, D. Dobkin, W. Sweldens and P. Schr\"oder,
Multiresolution mesh morphing, \emph{Proc. SIGGRAPH 99}, 1999.
\bibitem{LYS} C. Liu, Y. Yao and V. Srinivasan, Optimal process planning 
for laser 
forming of doubly curved shapes, \emph{J. Man. Sci. Eng.}, {bf 126} (2004) 1--9.
\bibitem{MH} J. Marsden and  T. Hughes, \emph{Mathematical Foundations of
Elasticity}, Prentice-Hall, Inc., 1983.
\bibitem{M} J. Moser,
On the volume elements on a manifold.
\emph{Trans. Amer. Math. Soc.}  {\bf 120} (1965) 286--294.
\bibitem{Hicks} N. Hicks, \emph{Notes on Differential Geometry},
Princeton NJ, D. Van Nostrand Company, Inc, 1965.
\bibitem{HLZ} Shi-Min Hu, Chen-Feng Li and Hui Zhang. 
\emph{Actual Morphing: A Physics-Based Approach to Blending}. ACM Symposium
on Solid Modeling and Applications (2004).
\bibitem{KN} S. Kobayashi and K. Nomizu, 
\emph{Foundations of differential geometry}, Interscience Publishers, 1963.
\bibitem{Wehrheim} K. Wehrheim,
\emph{Uhlenbeck Compactness}, European Mathematical Society, 2004.
\bibitem{W}  G. Wolberg, \emph{Digital Image Warping}, Los Alamitos, IEEE Computer Society Press, 1990.
\bibitem{W1} G. Wolberg, Image Morphing: A Survey, \emph{ Visual Computer}, {\bf  14} (1998),  360--372.
\bibitem{YPM} G. Yu, M. Patrikalakis and T. Maekawa,
Optimal development of doubly curved surfaces,
\emph{Computer Aided Geometic Design} {\bf 17} (2000) 545--577.
\end{thebibliography}
\end{document}